# Measuring Diversity
## A review and an empirical analysis

FRANCISCO PARREÑO
Escuela Superior de Ingeniería Informática, Universidad de Castilla La Mancha, Spain
Francisco.Parreno@uclm.es

RAMÓN ÁLVAREZ-VALDÉS
Departamento de Estadística e Investigación Operativa, Universidad de Valencia, Spain
Ramon.Alvarez@uv.es

RAFAEL MARTÍ
Departamento de Estadística e Investigación Operativa, Universidad de Valencia, Spain
Rafael.Marti@uv.es

**ABSTRACT**

Maximum diversity problems arise in many practical settings from facility location to social networks, and constitute an important class of NP-hard problems in combinatorial optimization. There has been a growing interest in these problems in recent years, and different mathematical programming models have been proposed to capture the notion of diversity. They basically consist of selecting a subset of elements of a given set in such a way that a measure based on their pairwise distances is maximized to achieve dispersion or representativeness. In this paper, we perform an exhaustive comparison of four mathematical models to achieve diversity over the public domain library MDPLIB, studying the structure of the solutions obtained with each of them. We extend this library by including new Euclidean instances which permit to analyze the geometrical distribution of the solutions. Our study concludes which models are better suited for dispersion and which ones for representativeness in terms of the structure of their solutions, as well as which instances are difficult to solve. We also identify in our conclusions one of the models which is not recommended in any setting. We finalize by proposing two improvements, one related to the models and one to solving methods. The computational testing shows the value of the analysis and merit of our proposals.

## 1. Introduction

The challenge of maximizing the diversity of a collection of points arises in a variety of settings, and the growing interest of dealing with diversity resulted in an effort to study these problems in the last few years. The terms diversity, dispersion, and equity can be found in many optimization problems with a similar meaning, and broadly speaking, they consist of selecting a subset of elements of a given set maximizing a distance measure. As far as we know, the operations research literature has focused on optimizing one of these models, ignoring the others, and no comparison among them has been performed. In this paper we undertake to explore and compare the different models proposed to achieve diversity.

We believe that the term diversity is somehow ambiguous in the context of combinatorial optimization, and some problems seem to look for *dispersion* among the selected points, while others want to achieve some *representativeness*, where the selected points are representing a class in the given set. In their seminal paper on diversity problems, Glover and Kuo (1998) say that "diversity is a rather nebulous term with overtones with a vaguely statistical nature". The following classic example mentioned in many diversity papers and books (Duarte et al., 2018) illustrates this point:

> *"Consider a group of university students from which we want to select five to form a committee. For each pair of students, we can compute a distance value based on their particular attributes. These include personal characteristics such as age or gender. A pair of students with similar attributes receives a low distance value, while students with different attributes have a large distance value. To create a diverse committee, the five students have to be selected in a way that their ten pairwise distances are relatively large."*

It is clear that in the example above, by a "diverse committee" we refer to a set of students that represents well the entire group in that university. In other words, it seems that the objective is that the selected students contain most of the attributes in the group, rather than their dispersion in the group. On other applications, such as telecommunications or transportation, we typically look for diverse routes or nodes in the sense of geographically disperse. We therefore consider these two classes of applications of diversity models: *representativeness* and *dispersion*. However, the difference between these two types of applications are not reflected in the models used to solve them, and in spite of its practical significance the characteristics of the solutions obtained are ignored, and all previous studies just focus on the objective function value.

In this paper we analyze and compare the different diversity models proposed and we disclose which models induce dispersion and which ones representativeness, and what is the structure of the solutions obtained with them. In line with this, we have found that no study analyzes the large set of previously reported instances (MDPLIB) in terms of their potential difficulty and interest in terms of comparing heuristic methods (which is its main use). We therefore propose a unified study based on an empirical comparison of models and instances. Additionally, we extend this library by including some Euclidean based instances, which permit to geometrically represent the solutions.

The problem of maximizing diversity deals with selecting a subset of elements $M$ from a given set $V$ in such a way that the diversity among the elements is maximized (Glover et al., 1995). The two most popular diversity models in the literature consist of maximizing the sum of the diversity (*MaxSum*) and maximizing the minimum diversity (*MaxMin*) in $M$. As pointed out by Prokopyev et al. (2009), both of these objectives are based on measures of efficiency:

> "The maximum dispersion problem primarily focuses on operational efficiency of locating facilities according to distance, accessibility, impacts, etc. It also arises in various other contexts including maximally diverse/similar group selection."

The *MaxSum* and *MaxMin* literature includes extensive surveys (Ağca, Eksioglu, Ghosh, 2000; Erkut, Neuman, 1989; Kuo, Glover, Dhir, 1993), exact methods (Ağca et al., 2000; Ghosh, 1996; Pisinger, 2006; Martí et al., 2010), and heuristics (Ghosh, 1996; Hassin et al., 1997; Kincard, 1992; Ravi et al., 1994; Resende et al., 2010).

Prokopyev et al. (2009) introduced two other models in the context of diversity called equity models, which incorporate the concept of fairness among candidates. These models appear in different settings, such as urban public facility location, diverse/similar group selection, and sub-graph identification, in which one may address fair diversification or assimilation among members of a network. The *MaxMinSum* diversity problem maximizes the minimum aggregate dispersion among the chosen elements, while the Minimum Differential Dispersion model, *MinDiff*, minimizes extreme equity values of the selected elements. Martínez-Gavara et al. (2017) and Duarte et al. (2015) proposed heuristic methods for the MaxMinSum and MinDiff respectively. However, the practical significance of these two new models, and their incremental contribution with respect to the classic MaxSum and MaxMin models have not been studied or evaluated.

In this paper we first review in Section 2 the four models proposed to achieve diversity: MaxSum, MaxMin, MaxMinSum, and MinDiff. Then, Section 3 describes the instances collected in the MDPLIB benchmark (Martí and Duarte, 2010). Section 4 describes our extensive experimentation to analyze these models from a graphical point of view, comparing the structure of their respective solutions obtained with CPLEX. Section 5 is devoted to our numerical analysis to compare these four models, evaluating the instances in the library. As explained in our conclusions, some models are better suited for dispersion, while others reflect better the representativeness. Section 6 complements the analysis with two proposals, an improved method and a new model. The paper finishes with some conclusions, where we can highlight that the two recent models, MaxMinSum and MinDiff, have a small contribution with respect to the classic and simpler models, MaxSum and MaxMin. In particular, the MinDiff obtains solutions with a poor structure in terms of diversity and we recommend to avoid its use.

## 2. Mathematical models

Given a graph $G = (V, E)$, where $V$ is the set of $n$ nodes and $E$ is the set of edges, let $d_{ij}$ be the inter-element distance between any two elements $i$ and $j$, let $M \subseteq V$ be the set of $m$ selected elements. The two most popular dispersion problems in the literature consist in maximizing the sum of the diversity (*MaxSum*) and maximizing the minimum diversity (*MaxMin*) in $M$. In mathematical terms, the standard formulation is built upon binary variables, $x_i \in \{0,1\}$, indicating whether element $i$ is selected or not, for $i = 1, \dots, n$ as follows:

**Model 1. MaxSum**

$$\text{Maximize} \quad \sum_{i<j} d_{ij} x_i x_j$$

$$\text{subject to} \quad \sum_{i=1}^{n} x_i = m$$

$$x_i \in \{0,1\} \quad i = 1, \dots, n$$

The objective function of the MaxSum problem is trivially formulated as a quadratic binary problem. Kuo et al. (1993) use this formulation to show that the clique problem (which is known to be NP-complete) is reducible to the MaxSum. These authors also suggest the transformation of the quadratic binary model into a mixed integer program of the following form, where new variable $y_{ij}$ replaces the product $x_i x_j$ in the above formulation:

$$\text{Maximize} \quad \sum_{i<j} d_{ij} y_{ij}$$

$$\text{subject to} \quad \sum_{i=1}^{n} x_i = m$$

$$x_i + x_j - y_{ij} \leq 1 \quad 1 \leq i < j \leq n$$

$$-x_i + y_{ij} \leq 0 \quad 1 \leq i < j \leq n$$

$$-x_j + y_{ij} \leq 0 \quad 1 \leq i < j \leq n$$

$$y_{ij} \geq 0 \quad 1 \leq i < j \leq n$$

$$x_i \in \{0,1\} \quad i = 1, \dots, n$$

Martí et al. (2010) showed that the best formulation, in terms of the size of the instances solved with CPLEX, is the following one with the objective function decomposed in the sum of some w-values:

$$\text{Maximize} \quad \sum_{i<n} w_i$$

$$\text{subject to} \quad \sum_{i=1}^{n} x_i = m$$

$$-\bar{D}_i x_i + w_i \leq 0 \quad 1 \leq i \leq n-1$$

$$-\sum_{j=i+1}^{n} d_{ij} x_j + \bar{\bar{D}}_i (1 - x_i) + w_i \leq 0 \quad 1 \leq i \leq n-1$$

$$x_i \in \{0,1\} \quad i = 1, \dots, n$$

Where $\bar{D}_i = \sum_{j=i+1}^{n} \max(0, d_{ij})$ and $\bar{\bar{D}}_i = \sum_{j=i+1}^{n} \min(0, d_{ij})$.

An interesting variant of the MaxSum model, is the MaxMean model, in which the number of selected elements is not predefined. It consists in the maximization of the average dispersion. In mathemarical terms, the objective function is

$$\text{Maximize} \quad \frac{\sum_{i>j} d_{ij} x_i x_j}{\sum_i x_i},$$

and the constraints are the same than in the MaxSum model with the exception of the first one, which sets the number of elements to $m$, that is not present in this model.

**Model 2. MaxMin**

The MaxMin diversity problem, also known as $m$-dispersion problem (Chandrasekaran and Daughety, 1981), can be trivially formulated by simply considering its objective function over the set $M$:

$$\text{Maximize} \quad \min_{i,j \in M} d_{ij}$$

$$\text{subject to} \quad M \subseteq V$$

$$|M| = m$$

The MaxMin objective function can be formulated with binary variables in a similar way than the MaxSum above. Kuo et al. (1993) proposed to avoid the non-linearity in this model, coming from both the product of variables and the *min* function, as follows:

$$\text{Maximize} \quad w$$

$$\text{subject to} \quad \sum_{i=1}^{n} x_i = m$$

$$(C - d_{ij}) y_{ij} + w \leq C \quad 1 \leq i < j \leq n$$

$$x_i + x_j - y_{ij} \leq 1 \quad 1 \leq i < j \leq n$$

$$-x_i + y_{ij} \leq 0 \quad 1 \leq i < j \leq n$$

$$-x_j + y_{ij} \leq 0 \quad 1 \leq i < j \leq n$$

$$y_{ij} \geq 0 \quad 1 \leq i < j \leq n$$

$$x_i \in \{0,1\} \quad i = 1, \ldots, n$$

The value of $C$ in the second constraint is simply a very large constant value, larger than all the $d_{ij}$ values. In this way, if $y_{ij} = 1$, then this constraint is simplified to $C - d_{ij} + w \leq C$, thus to $w \leq d_{ij}$. Then, when maximizing $w$, it will match the minimum of the $d_{ij}$ values. Alternatively, if $y_{ij} = 0$, then this constraint is simplified to $w \leq C$, which does not limit $w$ given the definition of $C$.

Recently, Sayyady and Fathi (2016) solve an alternative model consecutively to obtain the optimal solution of the MinMax model. In particular, they consider the node packing problem, in which given a threshold value $l$, a graph $G(l)$ is defined with the set $V$ of $n$ nodes of graph $G = (V, E)$, and the set of edges $E(l) = \{(i,j) \in E: d_{ij} < l\}$. The node packing problem consists of finding a maximum cardinality subset of nodes so that no two nodes in this subset are adjacent to each other. It can be formulated in mathematical terms with binary variables, $x_i$, indicating if node $i$ is selected as:

$$\text{Maximize} \quad \sum_{i=1}^{n} x_i$$

$$\text{subject to} \quad x_i + x_j \leq 1 \quad\quad 1 \leq i < j \leq n : d_{ij} < l$$

$$x_i \in \{0,1\} \quad i = 1, \ldots, n$$

The authors solve the node packing model above for different values of $l$ until they identify the optimal value of the MaxMin model. They propose an algorithm to strategically select the values of $l$ to solve a small number of node packing models.

Prokopyev et al. (2009) introduced several models to describe various aspects of diversity problems. In particular, these authors differentiate between efficiency and equity, and introduced mathematical programming formulations for the following equitable dispersion problems: MaxMinSum and MinDiff.

**Model 3. MaxMinSum**

The Maximum Minsum Dispersion Problem (*MaxMinSum*) consists of selecting a set $M \subseteq V$ of $m$ elements such that the smallest total dispersion associated with each selected element $i$ is maximized. The problem is formulated in Prokopyev et al. (2009) as follows:

$$\text{Maximize} \quad \{\min_{i: x_i=1} \sum_{j: j \neq i} d_{ij} x_j\}$$

$$\text{subject to} \quad \sum_{i=1}^{n} x_i = m$$

$$x_i \in \{0,1\} \quad i = 1, \ldots, n$$

The set of the $m$ selected elements is $M = \{i: x_i = 1\}$ and the objective function is based on measuring the total dispersion associated with each $i \in M$, denoted by $c(M, i)$. In particular, in this model $c(M, i) = \sum_{j \in M, j \neq i} d_{ij}$ and the objective is to maximize its minimum $c(M, i)$-value by a judicious selection of $M$.

The authors proposed an improvement over this model based on a big constant value $U^+$, which takes the value $U^+ = 1 + max_i \sum_{j \neq i} \max(0, d_{ij})$, and lower bounds $L_i = \sum_{j \neq i} \min(0, d_{ij})$ on the value of $\sum_{j \in M} d_{ij}$. We employ this tighter formulation in our empirical testing.

$$\text{Maximize} \quad s$$

$$\text{subject to} \quad s \leq \sum_{j: j \neq i} d_{ij} x_j - L_i(1 - x_i) + U^+(1 - x_i) \quad i = 1, \ldots, n$$

$$\sum_{i=1}^{n} x_i = m$$

$$x_i \in \{0,1\} \quad i = 1, \ldots, n$$

**Model 4. MinDiff**

Equity problems are an actual concern in the context of facility location, where the fairness among candidate facility locations is as relevant as the dispersion of the selected locations (Teitz, 1968). Given a set of elements, the main objective of these problems is to find a subset of those elements that minimizes a similarity measure. In this context, the MinDiff dispersion functions consider the extreme equity values of the selected elements. It can be formulated as:

$$\text{Minimize} \quad \max_{i \in M} \sum_{j: j \neq i} d_{ij} x_j - \min_{i \in M} \sum_{j: j \neq i} d_{ij} x_j$$

$$\text{subject to} \quad \sum_{i=1}^{n} x_i = m$$

$$x_i \in \{0,1\} \quad i = 1, \dots, n$$

As in the MaxMinSum model, Prokopyev et al. (2009) proposed an improved mixed linear 0-1 formulation based on bounds. Let $L_i$ and $U_i$ be lower and upper bounds on the value of $\sum_{j \in M} d_{ij}$ respectively. Let $U^+$ be an upper bound on the $U_i$ values, and let $L^-$ be a lower bound on the $L_i$ values.

$$\text{Minimize} \quad t$$

$$\text{subject to} \quad t \geq r - s \quad i = 1, \dots, n$$

$$r \geq \sum_{j: j \neq i} d_{ij} x_j - U_i(1 - x_i) + L^-(1 - x_i) \quad i = 1, \dots, n$$

$$s \leq \sum_{j: j \neq i} d_{ij} x_j - L_i(1 - x_i) + U^+(1 - x_i) \quad i = 1, \dots, n$$

$$\sum_{i=1}^{n} x_i = m$$

$$x_i \in \{0,1\} \quad i = 1, \dots, n$$

## 3. Instance Sets

Martí and Duarte (2010) compiled a comprehensive set of benchmark instances representative of the collections used for computational experiments in the MDP. They called this benchmark MDPLIB, which collects a total of 315 instances, and it is available at http://grafo.etsii.urjc.es/optsicom/. It basically contains three sets previously used: SOM, GKD, and MDG. A brief description of each set follows.

**SOM**: This data set consists of 70 matrices with integer random numbers between 0 and 9 generated from an integer uniform distribution.

- SOM-a: These 50 instances were generated by Martí et al. (2010) with a generator developed by Silva et al. (2004). The instance sizes are such that for n = 25, m = 2 and 7; for n = 50, m = 5 and 15; for n = 100, m = 10 and 30; for n = 125, m = 12 and 37; and for n = 150, m = 15 and 45.

- SOM-b: These 20 instances were generated by Silva et al. (2004) and used in most of the previous papers (see for example Aringhieri et al. 2008). The instance sizes are such that for n = 100, m = 10, 20, 30 and 40; for n = 200, m = 20, 40, 60 and 80; for n = 300, m = 30, 60, 90 and 120; for n = 400, m = 40, 80, 120, and 160; and for n = 500, m = 50, 100, 150 and 200.

**GKD:** This data set consists of 145 matrices for which the values were calculated as the Euclidean distances from randomly generated points with coordinates in the 0 to 10 range.

- GKD-a: Glover et al. (1998) introduced these 75 instances in which the number of coordinates for each point is generated randomly in the 2 to 21 range. The instance sizes are such that for n = 10, m = 2, 3, 4, 6 and 8; for n = 15, m = 3, 4, 6, 9 and 12; and for n = 30, m = 6, 9, 12, 18 and 24.

- GKD-b: Martí et al. (2010) generated these 50 matrices for which the number of coordinates for each point is generated randomly in the 2 to 21 range and the instance sizes are such that for n = 25, m = 2 and 7; for n = 50, m = 5 and 15; for n = 100, m = 10 and 30; for n = 125, m = 12 and 37; and for n = 150, m = 15 and 45.

- GKD-c: Duarte and Martí (2007) generated these 20 matrices with 10 coordinates for each point and n = 500 and m = 50.

**MDG**: This data set consists of 100 matrices with real numbers randomly selected between 0 and 10 from a uniform distribution.

- MDG-a: Duarte and Martí (2007) generated these 40 matrices, 20 of them with n = 500 and m = 50 and the other 20 with n = 2000 and m = 200. These instances were used in Palubeckis (2007).

- MDG-b: This data set consists of 40 matrices generated by Duarte and Martí (2007). 20 of them have n = 500 and m = 50, and the other 20 have n = 2000 and m = 200. These instances were used in Gallego et al. (2009) and Palubeckis (2007).

- MDG-c: Martí et al. (2013) proposed this data set with 20 matrices with n = 3000 and m = 300, 400, 500 and 600. These are the largest instances reported in our computational study. They are similar to those used in Palubeckis (2007).

The main objective of this paper is to evaluate the measures and models for diversity and not to compare the performance of heuristic solving methods, which has been extensively done in previous works. Our approach basically consists of solving the instances with CPLEX to optimality, size permitting. We therefore exclude from our study those very large instances ($n > 500$) that cannot be solved in this way. In particular, most of the MDG instances. To obtain results in this set with CPLEX, we created some new instances of a smaller size based on the existing ones. In particular, MDG-a2 contains 20 instances with $n = 100$ and $m = 10$ obtained by selecting the first 100 elements in the 20 instances of MDG-a of size 500. Similarly, MDG-b2 contains 20 instances with $n = 100$ and $m = 10$ obtained by selecting the first 100 elements in the 20 instances of MDG-b of size 500.

An important aspect that has been overlooked in most previous studies is the solution structure and the graphical representation of the elements selected in the solution. Empirical analyses in previous papers usually limit themselves to the comparison of some statistics based on the objective function and running times. We complement them with the present study, in which we compare not only the objective functions, but also the disposition of the points identified in the solution of each model.

Although the instances in the GKD set contain Euclidean distances computed from their points coordinates, these coordinates are not available, and the public records only contain the final distance matrices. Since we need the coordinates to represent the instances and their solutions, we have generated some additional Euclidean instances, to keep the coordinates and extend the MDPLIB by including them for future studies. For the sake of consistency with the existing library, we call them GKD-d.

- GKD-d: We generated 300 matrices for which the values were calculated as the Euclidean distances from randomly generated points with two coordinates in the 0 to 100 range. For each value of $n = 25, 50, 100, 250,$ and $500$, we generated 30 instances with $m = \lceil n/10 \rceil$ and 30 instances with $m = 2\lceil n/10 \rceil$.

An important characteristic of the matrices is the distribution of their values or even the number of different values that they may have. SOM instances represent an extreme case, artificially designed to pose a challenge for complex heuristic methods, only have integer values between 0 and 9. So each SOM matrix has a very large number of values repeated, and this is why it is a challenge for local search based methods that normally become lost in flat landscapes. On the contrary, GKD matrices have many different values. Being real numbers computed as a Euclidean distance with an accuracy of five decimal places, they have a large range of numbers. However, some of them are large and therefore they may also have repetitions. For example, a GKD-c instance is obtained from 10-coordinates points in the range [0,10], so if we have a point in (0,0,0,0,0,0,0,0,0,0) and another point in (10,10,10,10,10,10,10,10,10,10) then their distance is $\sqrt{10 \cdot 10^2} = 31.62$. That is clearly the maximum distance value in the matrix, and considering that a $500x500$ matrix has approximately 125,000 distance values, it is theoretically possible that all of them are different. However, we have empirically found when examining these instances that they have a normal distribution, and usually contain some repetitions, with a maximum value around 22.50.

Considering that the SOM instances have many repetitions, we wanted to generate instances with no repetitions in our GKD-d set. This is why we consider them with two [0,100] coordinates, because in this way, the maximum distance value is $\sqrt{2 \cdot 100^2} = 141.42$, and they hardly repeat the same value. Our empirical study will confirm the different performance of the models in the different types of instances

## 4. Graphical analysis of diversity models

In this section we solve the different models described in Section 3 on the GKD-d instances to represent their solutions, and analyze their structure. We start with the Min-Diff model since it presents the most elaborated mathematical objective function among the models considered, and it is especially difficult to anticipate the properties of its solutions. Then, we analyze in Section 4.2 the two models based on the sum function together, namely the MaxSum and the MaxMinSum. Finally, in Section 4.3 we review the MaxMin model.

## 4.1 The MinDiff model

Figure 1 shows four instances of the GKD-d set with $m$=25 and $m = 3$. In particular, they show the 25 points in the set (depicted with small black circles), and the 3 points selected in the optimal solution of the Min-Diff model (depicted with larger red circles). The objective function value in the solutions of the two bottom instances are very similar. Specifically, the instance represented in the bottom-right diagram has a value of 1.057, while the bottom-left instance has a value of 1.061. However, the disposition of the 3 selected points is very different in each one. In the bottom-left instance, the three points are located in the upper part of the diagram and are relatively close to each other. In the bottom-right one, they are in the outer part of the diagram, close to its border lines, and relatively far from each other.

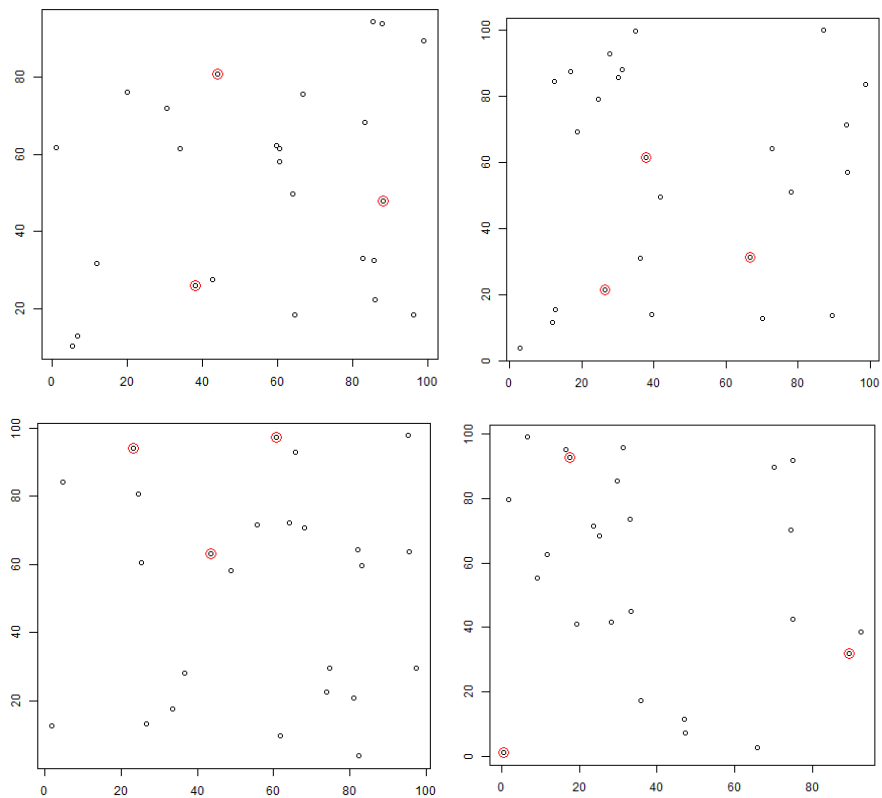

**Figure 1**. Optimal MinDiff solutions of four instances with $n = 25$ and $m = 3$.

In the objective function of the MinDiff problem, we want to minimize the difference between the maximum and the minimum of the distance values $d_i = \sum_{j \in M: j \neq i} d_{ij} x_j$ computed for each selected vertex ($i \in M$). The minimum objective function value of 0 can be therefore achieved if the maximum and the minimum $d_i$ −values are the same. Note however, that this computation ignores if these two values, the maximum and the minimum, are relatively high or low, but it only considers its difference. This fact has important consequences in the relative position of the vertices in the optimal solution, since the model seeks for $m$ vertices in positions that provide very similar inter distance values, which will produce a small difference between their maximum and minimum values.

In geometrical terms, we could say that if we are looking for $m = 3$ points in the plane to optimize the MinDiff model, their best possible location would be the vertices of an equilateral triangle, regardless if they are close or far to each other. If we join with lines the three points in any of the solutions in Figure 1, we will obtain a polyhedron similar to an equilateral triangle. In line with this, if we solve now the MinDiff model to select $m = 5$ points, we will see that the solution seeks for 5 elements with equal or at least very similar inter-distance values, which ideally would be obtained in the corner points of a regular pentagon. Figure 2 shows the optimal MinDiff solution of an instance with $m$=25 and $m = 5$, where we can confirm this pattern.

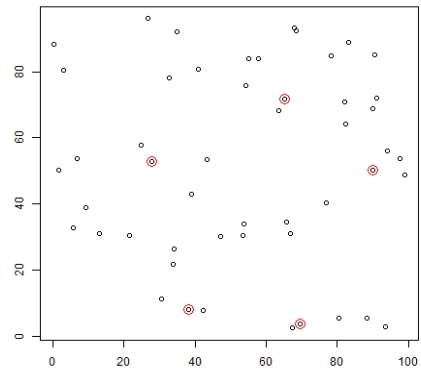

**Figure 2**. Optimal MinDiff solution for $n = 25$ and $m = 5$.

In the four examples depicted in Figure 1, the best MinDiff solution is in the top-right part of the figure. It has a value of 0.154, which compares favorably with the 0.361 of the top-left figure and with the other two examples in the bottom (1.057 for the bottom-right and 1.061 for the bottom-left). However, in spite of its larger objective function value (in a minimization problem), the selected points in the bottom-right figure are more disperse than those in the best solution. This is easily explained by the fact described above, this objective function only seeks for inter-distance equality among the selected points, and ignores how large or small they are. From this analysis, we can conclude that this model does not induce diversity or dispersion.

To complement this analysis, we represent in a single diagram the 30 GKD-d instances with $n = 25$. Figure 3 shows the 30 optimal solutions (in red) obtained for the MinDiff model with $m = 3$.

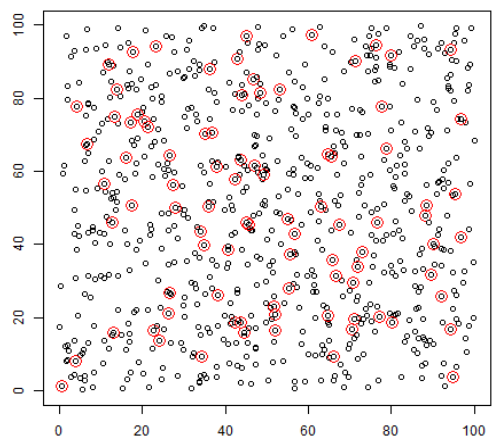

**Figure 3**. Optimal MinDiff solutions of 30 instances ($n = 25, m = 3$).

It is clear from Figure 3 that the solutions do not follow any specific pattern and they seem scattered all over the plane. To analyse them, we compute for each one the average inter-distances (i.e., the average of the pairwise distances of the selected elements). This value ranges in the 30 solutions from 10.725 to 94.158, which confirms our analysis that this model ignores if the inter-distances are small or large, but only focus on their similarity. This figure also illustrates that the MinDiff model does not induce representativeness in the sense that the selected points do not seem to be locate in the "center" of the set of points. To quantify it, we evaluate for each solution the average of the distance between its selected points and the non-selected points. This value ranges in the 30 solutions from 1.463 to 129.341, confirming the absence of a pattern in the location of the selected points both with respect to the other selected points, and to the non-selected ones.

It must be noted that in the SOM instances of the MDPLIB, we can find pairs of points at a distance of 0. Then, if we are solving the MinDiff, say for example for $m = 3$, and we found three points, in which the three pairwise distances are 0, the objective function would be 0, and they would be optimal. We have empirically found that this happens in many occasions, and therefore this set cannot be used to evaluate the performance of the solving methods for the MinDiff. Authors have to acknowledge that this set was introduced in the context of the MaxSum model to pose a challenge for heuristics (Martí et al., 2013), and it cannot be used to test other models without a proper analysis. In particular, it is misleading to use it for testing the MinDiff.

## 4.2 The MaxSum and MaxMinSum models

In this section we undertake to examine in graphical terms the structure of the solutions obtained with the MaxSum and MaxMinSum models on the GKD-d instances. Considering that these two models are mainly based on the sum of the distances of the selected elements, we want to disclose if they present significant differences in the type of solutions produced.

Figure 4 shows the solutions obtained with the two models on an instance with $n = 25$ and $m = 3$. In particular, Figure 4(a) shows the three points selected with the MaxMinSum model, and Figure 4(b) the selection of the MaxSum model.

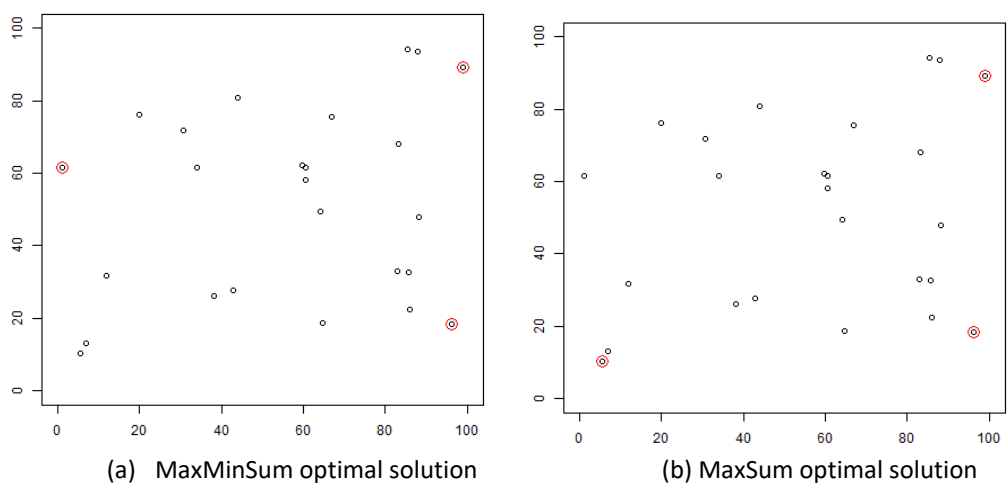

(a) MaxMinSum optimal solution　　　(b) MaxSum optimal solution

**Figure 4**. Optimal solutions of a GKD-d instance with $n = 25, m = 3$.

At a first glance, we can say that the two solutions shown in Figure 4 are similar, but not exactly the same. Specifically, both have two points in common out of the three ones selected. We can also analyze and compare their objective functions. Let $M^*$ be the optimal solution of the MaxSum model depicted in Figure 4(b), which has an objective function value in the MaxSum model of $z_{MS}(M^*) = 284.3$. Similarly, let $MM^*$ be the optimal solution of the MaxMinSum model depicted in Figure 4(a), which has an objective function value in the MaxMinSum of $z_{MMS}(MM^*) = 172.4$.

An interesting question is to compute how good is the solution obtained when solving one model with respect to the objective of the other model. If we compute the objective function value of $M^*$ in the MaxMinSum model, we obtain $z_{MMS}(M^*) = 162.0$, which has a relative percentage deviation with respect to the optimal solution $MM^*$ of this model of 6.03%. Similarly, we can compute the objective function value of $MM^*$ in the MaxSum model, obtaining $z_{MS}(MM^*) = 278.88$, which has a relative percentage deviation with respect to the optimal solution $M^*$ of this model of 1.91%. These results seem to indicate that the optimal solution obtained with one model scores relatively well in the other model, presenting a small deviation with respect to its optimum.

Figure 5 shows the same diagrams than Figure 4 but on a larger instance. In particular they show the MaxMinSum optimal solution, Figure 4(a), and MaxSum optimal solution, Figure 4(b) of a GKD-d instance with $n = 50$ and $m = 10$.

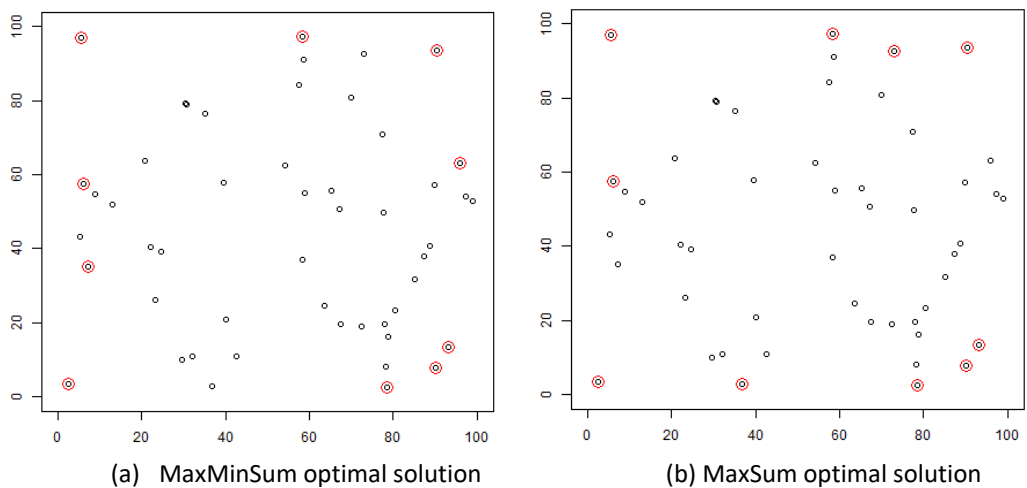

(a) MaxMinSum optimal solution     (b) MaxSum optimal solution

**Figure 5**. Optimal solutions of a GKD-d instance with $n = 50, m = 10$.

If we compare the two solutions in Figure 5, we reach the same conclusion than in Figure 4; both models provide similar solutions, but not exactly the same. In this case, their solutions share 8 points out of the 10 selected ones. If we also compare their objective values, and how well do they score in the "other model", we also confirm our previous impression: they work very well in the other model. Specifically, the optimal solution $M^*$ of the MaxSum model with $z_{MS}(M^*) = 3429.88$, presents a value in the MaxMinSum model of $z_{MMS}(M^*) = 638.27$, which represents a percentage deviation of its optimal value of 0.2%. Similarly, the optimal solution $MM^*$ of the MaxMinSum model with $z_{MMS}(MM^*) = 639.92$, presents a value in the MaxSum model of $z_{MS}(M^*) = 3426.25$, which represents a percentage deviation of its optimal value of 0.1%. Note that these percentage deviations are smaller than those obtained with the solutions in Figure 5. We performed additional experiments and confirm that the larger the size of the instances, the smaller these deviations are.

A geometrical aspect that we found relevant is that the 10 points in the two solutions of Figure 5 are located close to the border of the diagram, avoiding the central region of the plane. Figure 4 seems to confirm this disposition of the vertices; however, with just 3 selected points, it is difficult to assess. We therefore represent another instance of a large size to further investigate this point. Figure 6 shows the same diagrams than the previous figures on an instance with with $n = 100$ and $m = 20$. We did not obtain the optimal solution with CPLEX in the time limit of 1 hour of computer time, and the solutions represented in this figure are the best ones found within that time limit.

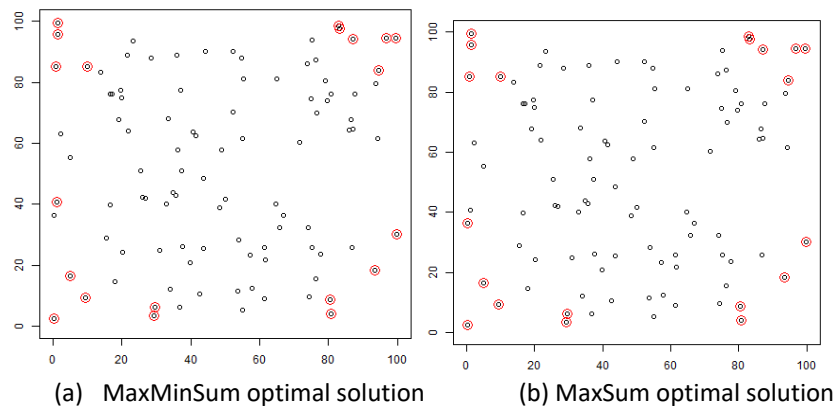

(a)   MaxMinSum optimal solution           (b) MaxSum optimal solution

**Figure 6**. Optimal solutions of a GKD-d instance with $n = 100, m = 20$.

The geometrical configuration of the selected points in Figure 6 confirms our conjecture that the MaxSum and MaxMinSum models obtain solutions close to the borders, and with no points in the central region. This clearly induces diversity or dispersion in the set of selected points, but the fact that the models do not select any solution in the central region can be also interpreted as a lack of representativeness. On the other hand, note that these two solutions only differ in the selection of one point, thus confirming as well that they are very similar.

Figures 5 and 6 also show that some of the pairs of the selected points are very close. These two models compute the overall sum of distances between the pairs of selected points. When we maximize that sum, we are looking for a selection in which these pairs are far apart. However, we can observe in the solutions represented that although most of them are far apart, some of them are very close. As a matter of fact, it is well-known in optimization that models based on a sum of items seek for a balance in the solution, in which the items with larger contribution compensate those with a smaller one, obtaining in that way the overall optimal configuration. However, in our problem in which we are looking for dispersion, the existence of some pairs in which the selected points are very close to each other, is an important drawback of these two models.

For the sake of brevity, we do not reproduce here all the drawings obtained with these two models on the different instances, but all of them show this type of disposition. We summarize them in Figure 7 in which we represent together all the solutions of each model for the set with $n = 25, m = 3$.

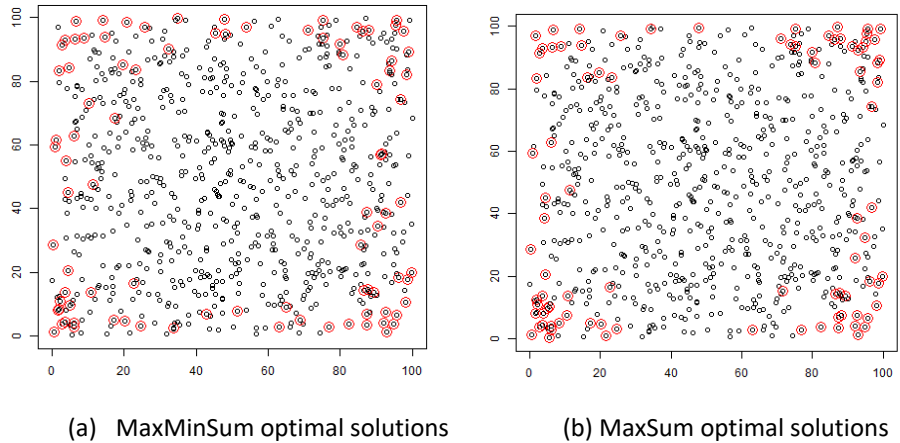

(a) MaxMinSum optimal solutions      (b) MaxSum optimal solutions

**Figure 7**. Optimal solutions of GKD-d instances with $n = 25, m = 3$.

Figure 7(a) shows the solutions (selections of 3 points) of the MaxMinSum model on the 30 GKD-d instances with 25 points. Figure 7(b) shows the 30 solutions of the MaxSum model on the same instances. The pairwise distance values between the 30 MaxMinSum optimal solutions ranges from 1.50 to 131.08, while the distances between the 30 MaxSum optimal solutions ranges from 1.50 to 129.43. Comparing both figures and their distances, it is clear that the solutions obtained with these two models are very similar.

We now repeat the experiment (i.e., represent in the same diagram the optimal solutions of many instances) on a larger instances set. Figure 8 shows the solutions of both models on an the subset of instances with $n = 50, m = 10$. The two diagrams in Figure 8 are in line with the previous ones and confirm the comments above on the relative location of the selected points and the similarities of both models.

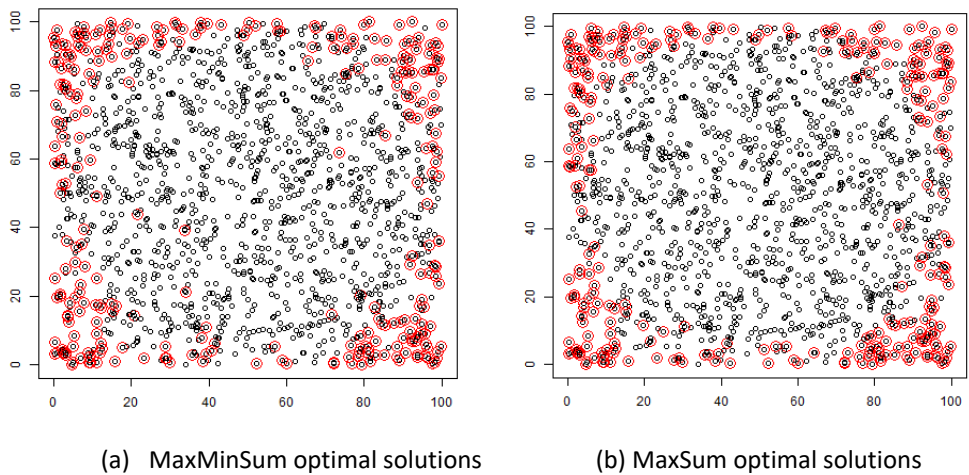

(a) MaxMinSum optimal solutions      (b) MaxSum optimal solutions

**Figure 8**. Optimal solutions of GKD-d instances with $n = 50, m = 10$.

## 4.3 The MaxMin model

In line with the study in the previous sections, we represent the optimal solution of a GKD-d instance. We start with the small one already studied, with $n = 25, m = 3$, so we can compare the MaxMin model solutions with the solutions obtained with the previous models. Figure 9 shows the MaxMin optimal solution.

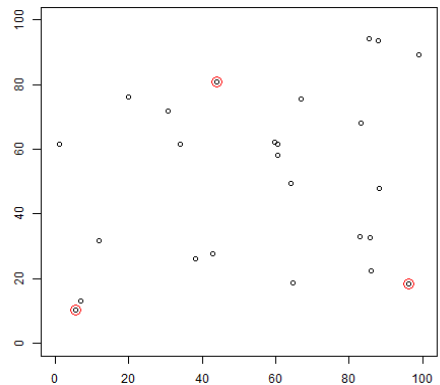

**Figure 9**. MaxMin optimal solution with $n = 25, m = 3$.

The optimal solution shown in Figure 9 is very similar to the ones obtained in Figure 4 for the MaxSum and MaxMinSum models, with three points selected in the outer part of the diagram. We consider now larger instances as we did above for the previous models.

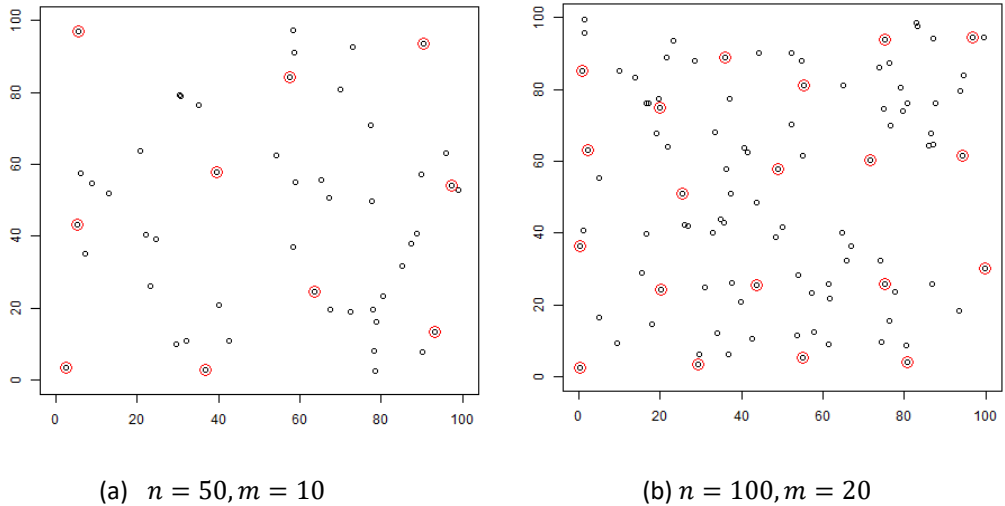

(a) $n = 50, m = 10$        (b) $n = 100, m = 20$

**Figure 10**. MaxMin optimal solutions of GKD-d instances.

The two solutions depicted in Figure 10 show a very different pattern compared with the previous solutions showed above. It is clear that this model does not avoid to select points in the central region of the plane, but on the contrary, it obtains equidistant points all over the region where the initial set of points rely. One could argue that it is not clear if it induces dispersion, since for each selected point we can find another one relatively close it. On the other hand, we may say that it induces representativeness in the sense that all the regions are represented in the selection of points. To confirm this point, we consider now the experiment in which we accumulate in the same diagram many optimal solutions.

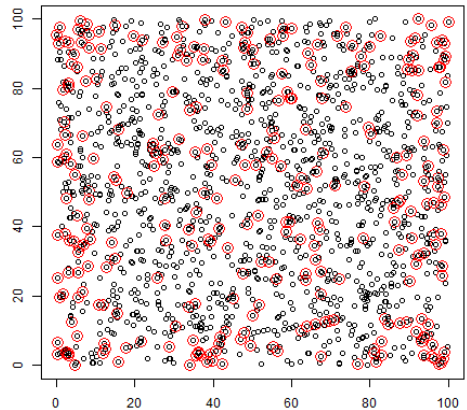

**Figure 11**. Optimal solutions of GKD-d instances with $n = 50, m = 10$.

Figure 11 shows the optimal solutions obtained in the 30 instances with 50 elements, where we select 10 of them. This graphic confirms that the MaxMin model does not avoid the selection of points in the central region, and it is meant to pick up equidistant points all over the region.

## 5. Numerical analysis of diversity models

In this section we solve the different models described in Section 2 on the entire benchmark of 675 instances collected in Section 3. We consider the most efficient formulations described for each model. We exclude from this analysis the original sets MDG in the MDPLIB because they are too large to be solved with CPLEX, and we include our two new sets MDG-a2 and MDG-b2 with 20 instances, each one of size $n = 100$. We employ CPLEX version 12.8.0.0 on a Red Hat 4.8.5-16Linux OS machine with 4 CPU, 8 Threads, 2.40GHz, and 4GB of RAM. The algorithms required to run the models in CPLEX were implemented in C++ on Visual Studio 2017 and compiled with g++.

| Instances | | MaxMin | | MaxMinSum | | MaxSum | | MinDiff | |
|---|---|---|---|---|---|---|---|---|---|
| Type | Num. | #Opt | Dev.(%) | #Opt | Dev.(%) | #Opt | Dev.(%) | #Opt | Dev.(%) |
| GKD-a | 75 | 75 | 0.0 | 75 | 0.0 | 75 | 0.0 | 75 | 0.0 |
| GKD-b | 50 | 50 | 0.0 | 49 | 0.03 | 15 | 79.3 | 20 | 59.9 |
| GKD-c | 20 | 0 | 24.0 | 0 | 6.8 | 0 | 584.5 | 0 | 100.0 |
| GKD-d | 300 | 300 | 0.0 | 187 | 6.4 | 100 | 136.0 | 144 | 47.2 |
| MDG-a2 | 20 | 20 | 0.0 | 20 | 0.0 | 0 | 114.9 | 0 | 100.0 |
| MDG-b2 | 20 | 20 | 0.0 | 20 | 0.0 | 0 | 129.9 | 0 | 100.0 |
| SOM-a | 50 | 50 | 0.0 | 30 | 5.3 | 15 | 76.2 | 15 | 70.0 |
| SOM-b | 20 | 11 | 40.0 | 1 | 37.3 | 0 | 158.9 | 0 | 100.0 |
| **Total** | **675** | **528** | **16.7** | **342** | **47.3** | **205** | **213.3** | **254** | **60.27** |

**Table 1.** Summary of results of the four models.

Table 1 shows, the name of each data set (column 1) and its number of instances (column 2). For each of the four diversity models, this table reports the number of optimal solutions obtained with CPLEX (#Opt), and the average percentage deviation of the best solution found with respect to the optimal solution. We limit the execution of CPLEX to a maximum of 3,600 seconds with each model on each instance.

Table 1 shows very interesting results when comparing the four models. The MaxMin model is able to solve most of the instances in our data set, while the MaxSum is the model that encounters more difficulties to solve them. Specifically, the MaxMin solves 528 out of the 675 instances in our study, which compares favorably with the 342 instances solved with the MaxMinSum, the 254 solved with the MinDiff, and the 205 solved with the MaxSum. This result is especially relevant when comparing the MaxSum and the MaxMinSum because we showed in the previous section that both obtain solutions of a similar structure, so in terms of obtaining the optimal solution with CPLEX, we would recommend to resort to the MaxMinSum. The remarkable performance of the MaxMin model is due to the transformation proposed by Sayyady and Fathi (2016) as a node packing problem. We empirically confirmed that previous MaxMin models exhibit a lower performance, in line with the other models.

Regarding problem size, most of the sets contain instances of medium size ($n = 100$), with the exception of the GKD-a that only contains very small instances ($n \leq 30$), and GKD-c with only large instances ($n = 500$). As expected, all the methods are able to solve to optimality all the instances in GKD-a, and none of them in GKD-c. Surprisingly, the MaxMin model is able to solve all the instances in the GKD-d set, which includes large instances of size $n = 500$.

It must be noted that, as can be seen in the Reference section, most of the papers devoted to diversity problems propose a heuristic and, in many cases, they test it on medium size instances. Our analysis reveals that if we are considering the MaxMin model, we can simply apply CPLEX in many cases to solve a medium size instance. As far as we know, the practical need of solving large size instances (with 500 elements or more) is not properly substantiated, since the most popular applications maximizing diversity, such as location (Erkut and Neuman, 1989), medical treatments (Kuo et al., 1993), or product design (Glover et al., 1998), are instantiated in relatively small size data. In this context, the common practice of comparing different metaheuristic technologies (Martí et al., 2013), although challenging, seems a rather artificial exercise only of theoretical relevance.

## 5. 1 Correlations

As described above, our point is that the MaxSum and MaxMinSum models provide very similar solutions in terms of their structure, which means that their optimal solutions share many elements. To confirm this point, we analyze now their objective function values both considering their deviations and correlations. To do that, we perform the test presented in Section 4.2 for an example, but now on the entire benchmark of 675 instances. Let $M^*$ be the optimal solution of the MaxSum model of a given instance, and $z_{MS}(M^*)$ its optimal value. Similarly, let $MM^*$ be the optimal solution of the MaxMinSum model, $z_{MMS}(MM^*)$ its optimal value, and $z_{MS}(MM^*)$ its value on the MaxSum model. Table 2 reports the average correlation between $z_{MS}(M^*)$ and $z_{MMS}(MM^*)$ for each set of instances. Additionally, we compute the percentage relative deviation of the MaxSum value of the optimal MaxMinSum solution. In mathematical terms:

$$Dev(MM^*) = 100 \cdot \frac{z_{MS}(M^*) - z_{MS}(MM^*)}{z_{MS}(M^*)}$$

Table 2 reports the average of the minimum, maximum and average $Dev$ values on each instance set. Note that when CPLEX is not able to obtain the optimal solution in one hour of computing time, we use the best solution found (i.e., a lower bound on the optimal value).

| Type | Num. | Correlation | Min. Dev. | Max. Dev. | Avg. Dev. |
|---|---|---|---|---|---|
| GKD-a | 75 | 1.00 | 0 | 2.29% | 0.41% |
| GKD-b | 50 | 1.00 | -0.02% | 1.70% | 0.41% |
| GKD-c | 20 | 0.99 | -0.02% | 0.33% | 0.13% |
| GKD-d | 300 | 1.00 | -0.09% | 14.91% | 1.00% |
| MDG-a2 | 20 | 0.00 | -0.56% | 5.22% | 2.16% |
| MDG-b2 | 20 | -0.07 | 1.83% | 5.44% | 3.44% |
| SOM-a | 50 | 1.00 | -5.75% | 5.88% | 0.17% |
| SOM-b | 20 | 1.00 | -9.33% | 2.13% | -1.53% |
| **Total** | **675** | 0.74 | -1.74% | 4.74% | 0.78% |

**Table 2.** MaxSum and MaxMinSum optimal solutions

Results in Table 2 confirm our point that MaxSum and MaxMinSum models are in general very similar considering that they obtain solutions of very related values. In fact, the average correlation between their optimal values is 1 in many instance sets, with the exception of the MDG sets, which present very low correlations (close to 0). On the other hand, the average relative deviation $Dev(MM^*)$ of the MaxMinSum optimal solution, $MM^*$, with respect to the MaxSum value is 0.78% on average over all the instances. In the MDGb2 instances this average value is 3.44% indicating the presence of some differences (although they can be still considered similar). Particularly interesting is the GKDd set, in which in spite of having an average of 1.0%, it presents a maximum value of 14.91%. This indicates that in some specific instances we can have optimal solutions in both models that are not similar, although they can be considered a kind of exception, with a general trend of being very similar.

### 5. 2 Graphical representations of solutions distances

We cannot represent in a direct way the elements (points in the space) and the solutions of randomly generated non-Euclidean instances, such as those in the SOM or MDG sets. We propose to represent and analyze the distance values of their solutions by depicting their distribution.

Given an instance in which we select $m = 5$ elements, a solution has 10 pairwise distance values between these 5 elements (note that it has a symmetric distance matrix, and therefore for each pair of selected elements we only compute one distance value). We perform an experiment in which we solve a model, say for instance the MaxMin, and collect the distance values in its optimal solution (for example the 10 pairwise distance values mentioned above). We consider a set of instances, for example the SOM set, solve all of them with a diversity model, and collect the pairwise distance values in all the optimal solutions obtained. In the case of the SOM instances, since the distances are integer values ranging from 0 to 9, we can represent them in a bar-chart with a bar for each number (category). Figure 12 shows the four bar charts corresponding to the four diversity models in our study for the SOM set of instances.

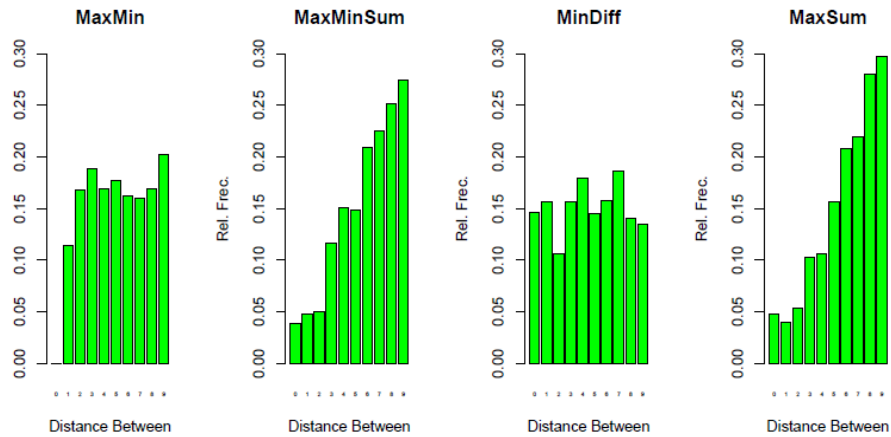

**Figure 12**. Distances in optimal solutions of SOM-a instances with $n = 50, m = 15$.

We may say in general terms, that a good solution to a diversity model should have large pairwise distance values and avoid the small ones. Therefore, we expect the histograms of the distances in the optimal solutions to be skewed left (i.e., with low frequencies in the first classes and larger ones in the upper classes, which is also called negatively skewed). Figure 12 shows that the histogram of the diversity models in the SOM instances are in line with this observation, although they present different patterns. The **MaxMin** is the only one of the four models that completely avoids small distances in its optimal solutions. We can see that the first class in the x-axis, which contains the smallest distances, is empty (i.e., with no distance values on it and thus with no bar). The relative frequency (depicted in the $y$ −axis) of the second class is around 0.1, and from the third to the 10$^{th}$ class we observe a very similar relative frequency (between 0.15 and 0.20). This indicates that the structure of the optimal solutions of this model completely avoids the points close to each other, but does not differentiate among points at a medium distance and points very far from each other, as shown in the flat profile of the histogram.

The bar chart of the **MaxSum** model presents a completely different profile than the one of the MaxMin model. In particular, the MaxSum exhibits a classic incremental distribution with a long tail, in which the relative frequency increases monotonically from the lower to the upper classes, presenting a negative skew. However, the smallest class is not completely avoided, as in the MaxMin model, and we can find very small distance values in its optimal solutions. In line with the results in the previous section, the **MaxMinSum** presents a very similar profile than the MaxSum, being very difficult to differentiate them.

In the analysis made in Section 4.1 on the **MinDiff** model, we conclude that this objective function only seeks for inter-distance equality among the selected points, and ignores how large or small they are. Figure 12 confirms this point since all the classes in the histogram have a similar frequency value, indicating that optimal solutions present all types of distances in a similar proportion, which clearly contradicts the seek for maximum inter-distances. We therefore ratify that this model does not induce diversity or dispersion.

We now turn our attention to the GKD instances. Since their distances take real numbers, we create here some histograms with the distance values of the optimal solutions in this set. To represent in a histogram all these distances together, we standardize them dividing these values by the maximum distance value in each corresponding instance. We consider the $x$-axis divided in 10 classes, where the first class contains the distances with a value lower than a 10% of the maximum value, the second class contains the distances between a 10% and a 20% of the maximum value, and so on. In this way, we analyze, across an entire set of instances, the distribution of the distances in the optimal solutions of a specific diversity model. Figure 13 shows the histograms of the inter-distance values of the optimal solutions on the GKD-b set. Results in this set are in line with those in the SOM set described above, although they present some differences.

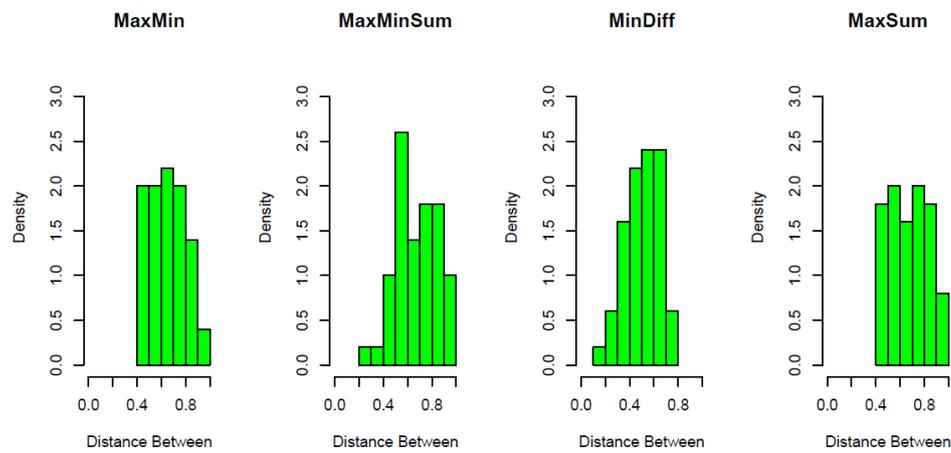

**Figure 13**. Distances in optimal solutions of GKD-b instances with $n = 50, m = 5$.

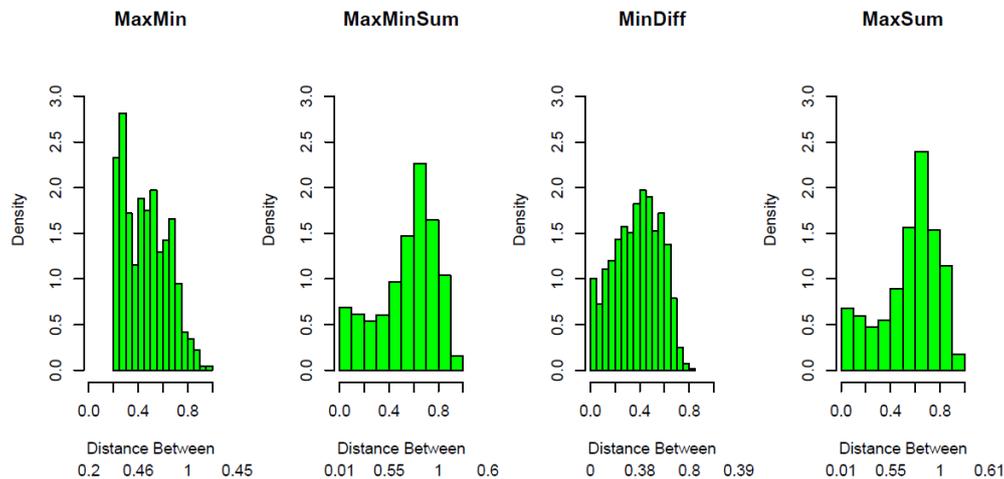

**Figure 14**. Distances in optimal solutions of GKD-d instances with $n = 50, m = 10$.

Figure 13 confirms that the MaxMin optimal solutions avoid the selection of low inter-distance values. The minimum value of the distances in these solutions is 0.44, while the minimum values in the MaxMinSum, MinDiff, and MaxSum optimal solutions are respectively 0.3, 0.17, and 0.42. Note the large value, 0.42, of the minimum distance in the MaxSum, which indicates a better performance of this model in the GKD than in the SOM instances, in terms of avoiding low distance values. However, the profile of the MaxMinSum and MaxSum solutions are not so clearly skewed as in the SOM instances, and it seems that they are not achieving in this set a large number of large inter-distance values. Finally, the MinDiff does not present a flat profile, and it seems more a Gaussian shape, indicating that it mostly selects central values.

Figure 14 represents the histograms of the four models on the GKD-d instances, in which we can see similar patterns than in the SOM instances. It is interesting to observe that the profiles of the histograms on the GKD-d are more similar to the SOM ones than to the GKD-b ones. In particular, the MaxMin model is the only one avoiding the selection of very low inter-distance values. Specifically, the minimum value in its optimal solutions is 0.2, which compares favorably with the 0.01 (MaxMinSum), 0 (MinDiff), and 0.01 (MaxSum). Profiles of the MaxSum and MaxMinSum are again left skewed, although their tail is not so long than in the SOM set, and the MinDiff exhibits a somehow mix profile between the one in the two sets above.

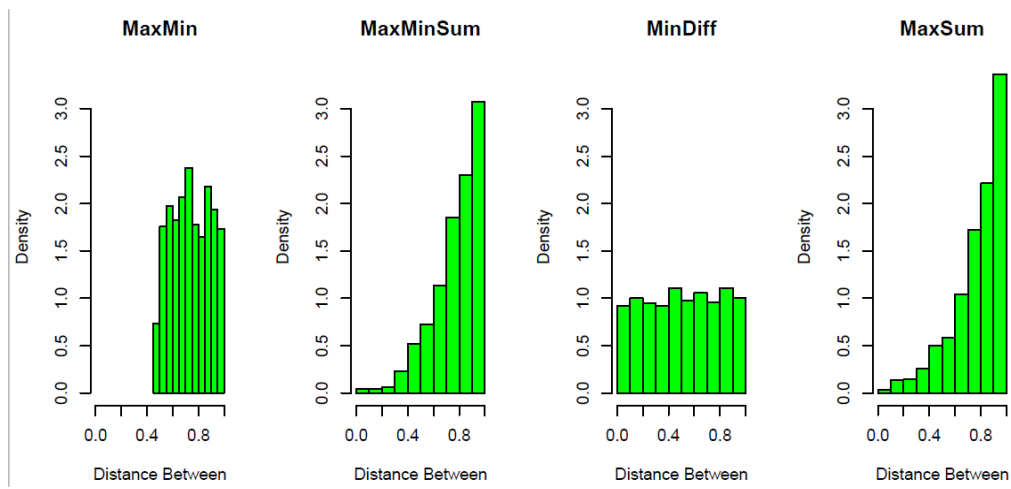

**Figure 15**. Distances in optimal solutions of MDG-a instances with $n = 100, m = 10$.

As commented above, there is a general trend in the profiles of the different models across the sets of instances, although there are differences. We conclude our analysis with the MDG-a that we found to be the most representative in which it shows a clear profile for each model. In particular, the MaxMin strongly avoids the small inter-distance values in its optimal solutions, with a minimum value of 0.46. The rest of the models present a minimum inter-distance value of 0. MaxMinSum and MaxSum are very similar and left skewed, with a mean (0.78) smaller than the median (0.83). Finally, the MinDiff distribution is completely flat, with a very similar frequency in each interval of the $x$-axis. We summarize now in Table 3 the main characteristics of each model.

| Model | Strengths | Drawbacks |
|---|---|---|
| MaxMin | Avoids very close elements. Optimal solutions can be found for medium size instances. | May select points either at a medium or at a large distance. |
| MaxSum | Favors the selection of points at a large distance | Permits very close elements |
| MaxMinSum | Similar to the MaxSum but permits to solve larger instances | Similar to the MaxSum but more complicated |
| MinDiff | Different structure than the other models, oriented to balance them. | Does not distinguish between close and apart elements, and does not induce diversity. |

**Table 3.** Models evaluation

## 6. Improvements

In this section we propose two improvements over the models and methods previously described. In particular, in the first subsection we propose a combined model of the MaxMinSum and MaxMin to overcome their respective limitations. In the second subsection we propose an enhancement on the Sayyady and Fathi (2016) exact method to solve the MaxMin.

### 6.1 A combined model

In the previous sections we clearly saw that the best models to achieve diversity are probably the MaxMin and the MaxSum (or similarly MaxMinSum), and both have strengths and limitations, as summarized in Table 3. We now try to put together the best characteristics of each model in a single one, overcoming their limitations at the same time. In particular, the MaxMin is good at avoiding very low inter-distance values but it does not discriminate between medium and large values. The MaxSum is somehow complementary because it is good at favoring the selection of large inter-distances, but it does not avoid the selection of very close elements.

**Model 5. Bi-level MaxSum**

Maximize $\sum_{i<j} d_{ij} x_i x_j$

subject to $\sum_{i=1}^{n} x_i = m$

$d_{ij} \geq d^* \quad \forall i, j \in M$

$x_i \in \{0,1\} \quad i = 1, \dots, n$

In our combined model, called Bi-level MaxSum, we propose to solve first the MaxMin problem to identify the maximum value $d^*$ of the minimum inter-distances in a solution with $m$ elements. Then, we solve the MaxSum model in which we add an extra constraint specifying that all the inter-distance values in the solution have to be larger than or equal to $d^*$. In this way, it avoids the small values that we observed in the standard MaxSum. On the other hand, we expect that it favors the selection of large inter-distance values, as in the standard MaxSum. In short, it presents the two best characteristics of the two models considered. It is clear however that we are solving a restricted version of the MaxSum model, and we thus expect that its objective function value deteriorates with respect to the optimum achieved with the standard formulation. Note that this can be considered a bi-level model according to Sinha et al. (2016), "Bi-level optimization is defined as a mathematical program, where an optimization problem contains another optimization problem as a constraint."

Additionally, note that the new model seeks for the best solution of the MaxSum in the restricted set of optimal solutions of the MaxMin. If the MaxMin does not had many alternate optimal solutions, the new model would have a small set of feasible solutions, and thus would be very constrained. We therefore perform an experiment to evaluate the number of alternate optimal solutions of the MaxMin as well as to quantify the deterioration in the MaxSum optimal value when comparing this restricted model to the original one.

Table 4 reports the alternate optimal solutions of the MaxMin model in the 80 instances with $n = 100$ and $m = 10$. In particular, for each set of instances with this size, it shows the number of instances (Num.), the average number of optimal solutions (Avg. Num.), and the maximum number of optimal solutions over the instances in the set (Max. Num.).

| Instances | Num. | Avg. Num. | Max. Num. |
|---|---|---|---|
| GKD-b | 5 | 7.6 | 18 |
| GKD-d | 30 | 317.9 | 8455 |
| MDG-a2 | 20 | 1.6 | 4 |
| MDG-b2 | 20 | 1.5 | 4 |
| SOM-a | 4 | 1040.7 | 2879 |
| SOM-b | 1 | 504.0 | 504 |
| **Total** | 80 | 178.8 | 8455 |

**Table 4.** Number of alternate optimal MaxMin solutions.

Table 4 shows that in some of the instance sets with $n = 100$ and $m = 10$, there are many alternate optimal MaxMin solutions, although in other sets this number is very small. The extreme cases are the GKDd set in which we can find an instance with 8455 optimal MaxMin solutions (with an average in this set of 317.9), and the MDGb2 in which all the instances have less than 5 alternate optimal solutions. Overall, the average number of alternate optima is 178.8, which supports our point of discriminating among them to identify the best one in terms of another objective function.

Considering that the MaxMinSum obtains solutions of a very similar structure than the MaxSum, and that it is easier to solve, as shown in Table 1, we consider the bi-level MaxMinSum in our next experiment. It can be formulated in a single mathematical model as:

**Model 6. Bi-level MaxMinSum**

$$\text{Maximize} \quad \{\min_{i:x_i=1} \sum_{j:j \neq i} d_{ij} x_j\}$$

$$\text{subject to} \quad \sum_{i=1}^{n} x_i = m$$

$$d_{ij} \geq d^* \quad \forall i, j \in M$$

$$d^* = \text{Max} \min_{i,j \in M'} d_{ij}$$

$$\text{subject to} \quad \sum_{i=1}^{n} y_i = m$$

$$y_i \in \{0,1\} \quad i = 1, \ldots, n$$

$$x_i \in \{0,1\} \quad i = 1, \ldots, n$$

Where $M'$ represents the solution of the lower level problem (MaxMin), defined by the $y$-variables, and $M$ represents the solution of the upper level problem (Constrained MaxMinSum) defined by the $x$-variables. Note that as a bi-level model, it is straightforward to solve, since the variables of the two problems are separated. We only have to solve the lower level model, MaxMin in this case, and then, with its optimal objective function value, solve the upper level model.

Table 5 reports the average percentage deviation (Avg. Dev.) of the MaxMinSum objective function of the bi-level model with respect to value of the original model. In other words, the relative deterioration of the objective function value when we reduce its feasible region to the optimal solutions of the MaxMin model.

| Instances | Num. | Avg. Dev. | Avg. Range. |
|---|---|---|---|
| GKD-b | 5 | 5.68 | 0.02 |
| GKD-d | 30 | 32.21 | 0.05 |
| MDG-a2 | 20 | 10.97 | 0.01 |
| MDG-b2 | 20 | 10.13 | 0.01 |
| SOM-a | 4 | 3.22 | 0.16 |
| SOM-b | 1 | 4.84 | 0.29 |
| **Total** | 80 | 17.93 | 0.04 |

**Table 5.** Bi-level MaxMinSum solutions

Table 5 shows that in some sets of instances, such as GKDb, the objective function value of the optimal solution presents a relative small reduction when we introduce the bi-level constraint in the MaxMinSum model. Specifically, the optimal value of the constrained model decreases in a 5.68% on average with respect to the optimal value of the original model. However, in other sets of instances, this reduction is more severe. In particular, the optimal value in the GKDd set experiences a 32.21% of reduction.

Table 5 also reports the results of another experiment related to the comparison of these two models, MaxMinSum (Model 3) and Bi-level MaxMinSum (Model 6). In particular, we collected all the optimal solutions of the MaxMin model and compute the value of the MaxMinSum function for all of them. The average range in Table 4 (last column) reports the average of the relative range of these values (computed as the maximum minus the minimum divided by the maximum value). The fact that these numbers are very low (lower than 0.3 in all cases) indicates that all these solutions present a similar value of the MaxMinSum function, and we can expect a marginal improvement when examining all of them in search for the optimum of the bi-level problem. This is why we consider an approximate method to heuristically solve the bi-level problem. We use a CPLEX function to collect a certain number of alternate optimal solutions of the MaxMin model, say 100,000, which is relatively efficient, even in large size instances (as shown in Table 1). Then, we limit the search of the bi-level problem to those solutions, and ignore the rest of them. It is clear that we don't have a proven optimal solution of the bi-level problem, but the result above indicates that we will obtain a solution with a value very close to optimality. Note that in the reported solutions we do not reach the 100,000 solutions, and therefore we obtain the optimal solution.

We now study the geometrical location of the selected points in the bi-level model. Since we are basically picking up one of the alternate optimal solutions of the MaxMin model, we may expect to obtain a solution with a similar structure to the one observed for this model. What it is not so predictable is the relative location of its selected points with respect to the ones selected in the MaxSum solution. Figure 16 shows the original MaxMin solution (a), the Bi-level MaxSum optimal solution (b), and the MaxSum optimal solution (c). a GKD-d instance with $n = 50, m = 5$

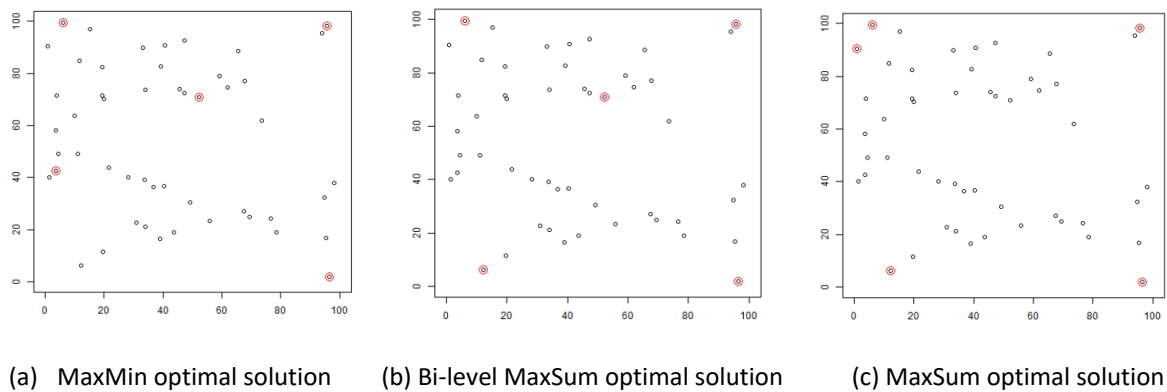

(a) MaxMin optimal solution    (b) Bi-level MaxSum optimal solution    (c) MaxSum optimal solution

**Figure 16.** Optimal solutions of GKD-d instance.

The MaxMin optimal solution depicted in Figure 16(a) shows the typical disposition represented in previous examples, in which the elements are scattered in the plane providing a disperse selection which includes the central region. A criticism of that selection however, would be the point in the left part of the diagram, around coordinates (5,40), instead of which we could easily select a better one in terms of global dispersion. As a matter of fact, the MaxSum value of that solution is 829.8, which is relatively low compared with the MaxSum optimal value of 942.8. The optimal MaxSum solution corresponding to that value is shown in Figure 6(c), and also has the typical disposition of that model, with the issue already mentioned of selecting two points very close (see the upper left corner of the square). Figure 16(b) clearly shows that the bi-level model provides an "in-between" solution considering the optimal solutions of the two original models. As expected, its solution is very similar to the MaxMin one, with only one different selected point, which in the new bi-level model we may say that is better located than in the original MaxMin in terms of favoring the total sum of distances. Instead of the point around coordinates (5,40), it selects the point around coordinates (15,5). This "exchange" does not change the MaxMin objective function, which is 51.4 in both models, but is able to increase the MaxSum value from 829.8 to 885.2 in the bi-level model. In short, we could say that the bi-level model tries to move points towards the outer part of the diagram, as the MaxSum model does, but avoiding to bring them together. We do not reproduce here other examples but we have found that the bi-level model is quite robust in terms of creating solutions with the best characteristics of the two models.

We would also like to highlight that the bi-level model can provide an efficient way to overcome the main limitation of the MinDiff model. As mentioned in previous sections, this objective function only seeks for inter-distance equality among the selected points, and ignores how large or small they are. Therefore, if we include a constraint specifying a minimum inter-distance value, the model would select a solution with balance inter-distances but necessarily large (according to that threshold). We do not include further experiments on this point for the sake of brevity, but it is clear that it is a point that deserves deep examination.

### 6.2 A more efficient algorithm for the MaxMin

As mentioned above, Sayyady and Fathi (2016) proposed an exact method for the MaxMin model based on solving several node packing problems. The authors introduced a parameter $l$, and defined an auxiliary graph $G = (V, E)$, over the same set of nodes than the original MaxMin problem, and set of edges $E(l) = \{(i,j) \in E: d_{ij} < l\}$. In this way, an optimal solution of the node packing problem in $G$ provides a set of points with minimum distance larger than or equal to $l$:

$v(l)$=Maximize $\quad \sum_{i=1}^{n} x_i$

subject to $\quad x_i + x_j \leq 1 \quad \forall (i,j) \in E(l)$

$\quad x_i \in \{0,1\} \quad i = 1, \dots, n$

Note however than in the MaxMin we specifically seek for a set of $m$ points, and the set obtained with the node packing has an arbitrary number of points, called $v(l)$. Sayyady and Fathi proposed to solve a sequence of node packing problems for different values of $l$ until we obtain a set of $v(l) = m$ points, which turns out to be the optimal solution of the MaxMin model.

The authors propose a systematic search in the interval $l \in [d_{min}, d_{max}]$, where $d_{min}$ and $d_{max}$ are the minimum and maximum values respectively among all the distances in the graph. They base the search on the property that if we have two $l-$values, $l_1$ and $l_2$ such that $l_1 < l_2$, $v(l_1) \geq m$ and $v(l_2) < m$ then, the optimal value of the MaxMin problem $z^*$ verifies $l_1 \leq z^* < l_2$. The method performs a binary search over the ordered set of different distances in the graph. To that end, the minimum distance between consecutive values is computed to divide the search into $2^q$ equal subintervals, performing at most $q$ steps (i.e., solving the node packing problem a maximum of $q$ times).

As disclosed our empirical comparison in Section 5, the effectiveness of the method by Sayyady and Fathi makes the MaxMin the most appealing diversity model, since it permits to solve relatively large problems to optimality. However, when implementing it, we found out two further improvements on this method that may speed up the process and therefore obtain better solutions when the time limit is reached.

Instead of solving the node packing problem, that returns a set of $v(l)$ points as shown above, we propose to include in the node packing model the standard size constraint of the diversity model (i.e., $\sum_{i=1}^{n} x_i = m$). We therefore do not consider the objective function in this model since it would be equal to $m$, and simply submit to CPLEX the modified model to check its feasibility, which is solved faster than the standard node packing. This is equivalent to answer the question:

¿Is there a solution $x_i \in \{0,1\}$ $i = 1, \dots, n$, such that $\sum_{i=1}^{n} x_i = m$, and $x_i + x_j \leq 1$ $\forall (i,j) \in E(l)$ ?

Additionally, instead of performing a binary search over the intervals (with a number larger than or equal to $2^q$), we perform a binary search directly over the set of distances $[d_{min}, d_{max}]$. We empirically found that the number of subintervals can be very large in medium size instances, and the original method does not consider the distribution of the distances, while we propose to divide the set of distances in a sequential fashion according to that distribution. In particular, we start by making $l$ to the median value in $[d_{min}, d_{max}]$, and solve the decision problem described above. If the answer is yes (i.e., if there is a solution with $m$ elements and value $l$), then we resort to the interval $[l, d_{max}]$; otherwise we consider $[d_{min}, l]$. We set now $l$ as the median value of the new interval and proceed in this way. Table 6 compare the performance of the original method and the new improved implementation. To simplify the study, we limit the MDG sets to the large examples with $n = 250$ and $m = 25$.

| Instances | | Original | | Improved | |
|---|---|---|---|---|---|
| Type | Num. | #Opt | Dev.(%) | #Opt | Dev.(%) |
| GKD-a | 75 | 75 | 0.0 | 75 | 0.0 |
| GKD-b | 50 | 50 | 0.0 | 50 | 0.0 |
| GKD-c | 20 | 0 | 24.3 | 0 | 8.9 |
| GKD-d | 300 | 300 | 0.0 | 300 | 0.0 |
| MDG-a | 20 | 0 | 71.20 | 0 | 19.8 |
| MDG-b | 20 | 0 | 76.86 | 0 | 18.7 |
| SOM-a | 50 | 50 | 0.0 | 50 | 0.0 |
| SOM-b | 20 | 11 | 40.0 | 17 | 15.0 |
| Total | 555 | 486 | 26.55 | 492 | 7.81 |

**Table 6.** MaxMin implementations

Table 6 clearly shows that the new implementation is able to improve the original one in those cases where it was not able to solve the instances. In particular, the original method presents an average percent deviation of 26.55, while the new one reduces it to 7.81 in our test sets. We generated 120 additional larger instances in the GKD set to perform a comprehensive comparison (note that the size of the instances is a critical factor here). In particular, we created 60 instances of 1,000 elements, and 60 of 2,000 elements. The original method was able to solve 42 of them, while the new one solves 76 out of the 120.

## Conclusions

This paper explores the four mathematical formulations studied in the recent literature to model diversity in the context of combinatorial optimization: MaxMin, MaxSum, MaxMinSum, and MinDiff. Our analysis is based on both numerical and geometrical representation. As far as we know, this is the first time that such a study analyzes the structure and properties of the solutions obtained with these models. Additionally, we propose a combined new model and improve a solution method.

The first conclusion of our study is that the **MaxSum** and **MaxMinSum** provide similar solutions, and considering the relatively large amount of research already done in the MaxSum model, it is not well justified the need of the recently introduced MaxMinSum one (especially because it is more complicated). In particular, the empirical analysis reveals that the optimal solution obtained with one model scores very well in the other model, presenting a small deviation with respect to its optimum (0.8% on average on the MDPLIB). Additionally, both models present an average correlation of 0.74, and in many cases it is larger than 0.9. Regarding the geometrical disposition of its solutions, we saw that they select points close to the borders of the space, and with no points in the central region. As an important drawback, these models may select very close elements (although they are usually just a few of them). In general terms we may say that they model what we usually understand as **dispersion** or diversity.

An important second group of conclusions is that the **MaxMin** model generates solutions with a very different structure than the MaxSum model. In particular, it obtains equidistant points all over the space, and it does not avoid to select points in the central part. On the other hand, if we want to obtain the optimal solution of an instance with CPLEX, the MaxMin model is able to solve most of the instances in our data set (528 out of the 675 instances), while the MaxSum is the model that encounters more difficulties to solve them (it only solves 205 to optimality). Based on the geometrical structure of its solutions, we can say that this mathematical formulation induces **representativeness**, more than dispersion, in their selection of points.

A third conclusion of our study is that the **MinDiff** only seeks for inter-distance equality among the selected points, and ignores how large or small these distances are. This is a very important observation since it disqualifies the term diversity or dispersion. The conventional wisdom usually interprets diversity as long distance and this model does not consider it. Although we may say that it balances the selection, achieving equity, we cannot find any reason to select balanced points at a very small distance. Our recommendation would be to **avoid the use of this model**, unless very well justified.

Our geometrical and numerical analysis are complemented with a visual inspection of the diagrams (bar-charts and histograms) of the distribution of distances in the optimal solutions of the four models. We conclude that the MaxMin is the only one of the four models that completely avoids small distances in its optimal solutions, and the MaxSum exhibits a classic incremental distribution in which the relative frequency increases monotonically from the lower to the upper classes.

Regarding the set of instances in the MDPLIB, we also conclude three interesting points. We first create a new set of Euclidean instances, GKD-d, and make them, including points coordinates, available for future geometrical studies. On the other hand, in the SOM instances, we can find pairs of points at a distance of 0. These instances are completely artificial and do not represent any real situation. They were created to compare heuristics and discriminate among them in the context of the MaxSum model. However, they may mislead empirical comparisons in other models, especially in the case of the MinDiff, in which we would recommend not to use them. Finally, we found the MDG-a the most representative set, in which we obtain with each model what is expected according to its definition, with clearly different solutions for each one.

We finish our study by proposing a new model, the Bi-level MaxSum, to capture the good properties of the best two models MaxMin and MaxSum. This model only consists of adding an additional constraint to the standard MaxSum. In particular, it avoids the small distance values that we observed in the standard MaxSum, while favoring the selection of large inter-distance values. Our empirical analysis shows that it is an improved diversity model. As a future line of research we believe that the inclusion of the additional constraint proposed in the Bi-level model could solve the MinDiff drawback identified in our study, thus obtaining a solid model to induce equity. We finished our paper by improving the iterative method to solve the MinMax model, which constitutes now the most effective method on a diversity problem.

## Acknowledgement

This research has been partially supported by the Spanish Ministry with grant ref. PGC2018-0953322-B-C21/MCIU/AEI/FEDER-UE.

## References

Ağca, S., B. Eksioglu, J. B. Ghosh; "Lagrangian solution of maximum dispersion problems". *Naval Research Logistics*, **47**: 97–114, 2000.

Duarte, A., M. Laguna, R. Martí; "MetaHeuristics for Business Analytics. A Decision Modeling Approach", EURO Advanced Tutorials on Operational Research, Springer, ISBN 978-3-319-68119-1, pp. 136, 2018.

Duarte, A., J. Sánchez-Oro, M. Resende, F. Glover, R. Martí; "GRASP with Exterior Path Relinking for Differential Dispersion Minimization", Information Sciences **296**, 46-60, 2015.

Erkut, E., S. Neuman; "Analytical models for locating undesirable facilities". *European Journal of Operational Research*, **40:** 275–291, 1989.

Erkut, E; "The discrete p-dispersion problem". *European Journal of Operational Research*, **46**:48–60, 1990.

Gallego, M., A. Duarte, M. Laguna, R. Martí; "Hybrid heuristics for the maximum diversity problem", Computational Optimization and Applications 44(3), 411-426, 2009.

Ghosh, J. B; "Computational aspects of the maximum diversity problem". *Operations Research Letters,* **19:** 175–181, 1996.

Glover, F., C. C. Kuo, K.S. Dhir; "A discrete optimization model for preserving biological diversity". *Applied Mathematical Modeling,* **19**: 696-701, 1995.

Glover, F., C. C. Kuo, K. S. Dhir; "Heuristic algorithms for the maximum diversity problem". *Journal of Information and Optimization Sciences,* **19**(1): 109-132, 1998.

Hassin, R., S. Rubinstein, A. Tamir; "Approximation algorithms for maximum dispersion". *Operations Research Letters*, **21:** 133–137, 1997.

Kincard, R. K; "Good solutions to discrete noxious location problems via metaheuristics". *Annals of Operations Research*, **40**: 265-281, 1992.

Kuo, C. C., F. Glover, K. S. Dhir; "Analyzing and modeling the maximum diversity problem by zero-one programming". *Decision Sciences*, **24:**1171–1185, 1993.

Martí, R., A. Duarte; The MDPLIB at Optsicom, http://grafo.etsii.urjc.es/optsicom/, 2010.

Martí, R., M. Gallego, A. Duarte, E. Pardo; "Heuristics and Metaheuristics for the maximum diversity problem", Journal of Heuristics **19** (4), 591-615, 2013.

Martí, R., F. Sandoya; "GRASP and path relinking for the equitable dispersion problem". *Computers and Operations Research,* **40**: 3091-3099, 2013.

Martí, R., M. Gallego, A. Duarte; "A Branch and Bound Algorithm for the Maximum Diversity Problem", European Journal of Operational Research **200**(1), 36-44, 2010.

Martínez-Gavara, A., V. Campos, M. Laguna, R. Martí; "Heuristic Solution Approaches for the Maximum MinSum Dispersion Problem, Journal of Global Optimization 67(3), 671-686, 2017.

Pisinger, D; "Upper bounds and exact algorithms for p-dispersion problems". *Computers and Operations Research*, **33**: 1380–1398, 2006.

Prokopyev, O. A., N. Kong, D. L. Martinez-Torres; "The equitable dispersion problem". *European Journal of Operation Research,* **197**: 59-67, 2009.

Ravi, S. S., D. J. Rosenkrantz, G.K. Tayi; "Heuristic and special case algorithms for dispersion problems". *Operations Research*, **42**: 299-310, 1994.

Resende, M. G. C., Martí, M. Gallego, A. Duarte; "GRASP and path relinking for the max–min diversity problem". *Computers and Operations Research*, **37**(3): 498-508, 2010.

Sayyady, F., Y. Fathi; "An integer programming approach for solving the p-dispersion problem". *European Journal of Operational Research*, 253: 216-225, 2016.

Sinha, A., P. Malo, and K. Deb. "Evolutionary algorithm for bilevel optimization using approximations of the lower level optimal solution mapping". European Journal of Operational Research 275(2), 2016